\documentclass[a4paper,10pt]{article}
\usepackage{stmaryrd}
\usepackage{amsfonts}
\usepackage{bbm}
\usepackage{amscd}
\usepackage{mathrsfs}
\usepackage{latexsym,amssymb,amsmath,amscd,amscd,amsthm,amsxtra,xypic}
\usepackage[dvips]{graphicx}
\usepackage[utf8]{inputenc}
\usepackage[T1]{fontenc}
\usepackage{lmodern}
\usepackage{amssymb}
\usepackage[all]{xy}
\usepackage{nicefrac,mathtools,enumitem}
\usepackage{microtype}

\newtheorem{thm}{Theorem}[section]
\newtheorem{lem}[thm]{Lemma}
\newtheorem{cor}[thm]{Corollary}
\newtheorem{pro}[thm]{Proposition}
\newtheorem{ex}[thm]{Example}
\newtheorem{rmk}[thm]{Remark}
\newtheorem{defi}[thm]{Definition}

\newcommand {\emptycomment}[1]{}

\newcommand{\be }{\begin{equation}}
\newcommand{\ee }{\end{equation}}

\newcommand{\pf}{\noindent{\bf Proof.}\ }

\newcommand{\huaV}{\mathcal{V}}

\newcommand{\D}{\mathfrak D\mathrm {er}}

\newcommand{\frkf}{\mathfrak f}
\newcommand{\g}{\mathfrak g}
\newcommand{\h}{\mathfrak h}
\newcommand{\frkg}{\mathfrak g}

\def\qed{\hfill ~\vrule height6pt width6pt depth0pt}

\newcommand{\half}{\frac{1}{2}}

\newcommand{\p}{\mathbbm{p}}
\newcommand{\br}[1]{   [ \cdot,    \cdot  ]   }
\newcommand{\id}{\mathbbm{i}}
\newcommand{\idd}{\mathrm{id}}

\newcommand{\dM}{\mathrm{d}}
\newcommand{\Z}{\mathrm{Z}}

\newcommand{\Hom}{\mathrm{Hom}}
\newcommand{\Der}{\mathrm{Der}}
\newcommand{\Pro}{\mathrm{Pr}}

\newcommand{\gl}{\mathfrak {gl}}

\newcommand{\ad}{\mathrm{ad}}

\newcommand{\sgn}{\mathrm{sgn}}
\newcommand{\Ksgn}{\mathrm{Ksgn}}

\begin{document}
\title{
{On non-abelian extensions of Leibniz algebras\thanks{This research is supported by NSFC (11471139) and NSF of Jilin Province (20140520054JH).}}
 }\vspace{2mm}
\author{  Jiefeng Liu$^1$, Yunhe Sheng$^2$  and Qi Wang$^2$\\
$^1$Department of Mathematics, Xinyang Normal University,\\ \vspace{2mm} Xinyang 464000, Henan, China\\
$^2$Department of Mathematics, Jilin University,\vspace{2mm}
 Changchun 130012, Jilin, China\\
  Email:~  jfliu13@mails.jlu.edu.cn,~ shengyh@jlu.edu.cn,~ wangq13@mails.jlu.edu.cn
}

\date{}

\maketitle

\begin{abstract}
In this paper, first we classify non-abelian extensions of Leibniz algebras by the second non-abelian cohomology. Then, we construct Leibniz 2-algebras using derivations of Leibniz algebras, and show that under a condition on the center, a non-abelian extension of   Leibniz algebras  can be described by a Leibniz 2-algebra morphism. At last, we give a description of non-abelian extensions in terms of Maurer-Cartan elements in a differential graded Lie algebra.
\end{abstract}
\section{Introduction}
Eilenberg and Maclane developed a theory of non-abelian extensions of abstract groups in
the 1940s \cite{Group extension}, leading to the low dimensional non-abelian group cohomology. Then there are a
lot of analogous results for Lie algebras
\cite{AMR,Hochschild,Nonabeliancohomology}.  Abelian extensions of Lie algebras can be
classified by the second   cohomology.  Non-abelian extensions of Lie algebras can be
described by some linear maps regarded as derivations of Lie algebras. This result is
generalized to the case of super Lie algebras in \cite{AMR2,Fialowski ext}, and to
the case of Lie algebroids in \cite{brahic}. In   \cite{nonabelin cohomology of Lie}, the author gave a description of non-abelian extensions of Lie algebras in terms of the  Deligne groupoid of a differential graded Lie algebra (DGLA for short). In \cite{shengzhu}, a non-abelian extension of the Lie algebra $\g$ by the Lie algebra $\h$ is explained by a Lie 2-algebra morphism from $\g$ to the strict Lie 2-algebra $\h\stackrel{\ad}\longrightarrow\Der(\h)$, where $\Der(\h)$ is the Lie algebra of derivations on $\h$. See \cite{baez:2algebras} for more information about a Lie 2-algebra, which is a categorification of a Lie algebra.


The notion of a Leibniz algebra was introduced by Loday \cite{Loday,Loday and Pirashvili}, which is a noncommutative generalization of a Lie algebra. Central extensions of Leibniz algebras are studied from different  aspects \cite{ccentral extension1,ccentral extension2,LCT of Lie}. Abelian extensions of Leibniz algebras is studied   in  \cite{cAbelian extension}.   Non-abelian cohomology and extensions by crossed modules of Leibniz algebras are studied in \cite{Casas,Gnedbaye, Khmaladze}. Associated to a Leibniz algebra, there is a differential graded Lie algebra structure on the graded vector space of the cochain complex, which plays an important role in studying cohomology and deformations  of Leibniz algebras, see \cite{dgLa1,dgLa2} for details.

 \emptycomment{Recently, Lie algebroids, Courant algebroids, and Poisson tensions ,etc. admit a algebraical description in terms of  differential graded Lie algebra. Frégier,Y also finds a good description of non-abelian extension of Lie algebras in this way. So the one of my present work  is to find a description of non-abelian extension of Leibniz algebras with the same method.}

 In this paper, we study non-abelian extensions of Leibniz algebras. First we introduce the second non-abelian cohomology $H^2(\g,\h)$ of a Leibniz algebra $\g$ with the coefficients in a Leibniz algebra $\h$, by which we classify non-abelian extensions of   $\g$ by   $\h$.  Then we construct two Leibniz 2-algebras using left derivations and right derivations of Leibniz algebras. One is a natural generalization of the Lie 2-algebra associated to derivations of a Lie algebra; The other can be used to describe non-abelian extensions of Leibniz algebras that satisfy a condition on centers. Whether there is a Leibniz 2-algebra associated to derivations of Leibniz algebras, that could be used to describe all the non-abelian extensions of Leibniz algebras, like the case of Lie algebras, remains to be an interesting problem. Finally, we generalize the result in \cite{nonabelin cohomology of Lie} to the case of Leibniz algebras. We construct a DGLA associated to Leibniz algebras $\g$ and $\h$, and describe non-abelian extensions of $\g$ by $\h$ by its Maurer-Cartan elements.

The paper is organized as follows. In Section $\ref{sec:pre}$, we study properties of left derivations $\Der^L(\g)$ and right derivations $\Der^R(\g)$ of a Leibniz algebra $\g$, and construct a semidirect product Leibniz algebra $\Der^L(\g)\oplus \Der^R(\g)$. Then we introduce the second  non-abelian cohomology of Leibniz algebras, by which we classify non-abelian extensions of Leibniz algebras (Theorem $\ref{classification}$). In Section $\ref{Leib20ext}$, first we recall some basic definitions regarding Leibniz $2$-algebras and morphisms between them. Then we construct Leibniz algebras $\Pi(\h)$ and $\Xi(\h)$ associated to a Leibniz algebra $\h$, which are subalgebras of $\Der^L(\h)\oplus \Der^R(\h)$. We go on constructing strict Leibniz 2-algebras $(\h,\Pi(\h),(\ad^L,\ad^R),l_2)$ and $(\h,\Xi(\h),(\ad^L,\ad^R),l_2)$. The former is a natural generalization of the Lie 2-algebra associated to derivations of a Lie algebra, and the latter can be used to describe some non-abelian extensions of $\g$ by $\h$. More precisely, we show that there is a one-to-one correspondence between non-abelian extensions $\hat{\g}$ of   $\g$ by $\h$ with $\Z(\h)=\Z(\hat{\g})\cap \h$ and morphisms from $\g$ to $(\h,\Xi(\h),(\ad^L,\ad^R),l_2)$   (Theorem $\ref{LB2toext1}$). In Section $\ref{dgla0ext}$, we construct a DGLA $(L,[\cdot,\cdot]_C,\overline{\partial})$, which is a sub-DGLA of $(C(\g\oplus\h),[\cdot,\cdot]_C,\overline{\partial})$, where $\g\oplus\h$ is the Leibniz algebra direct sum of $\g$ and $\h$. Then we show that  the equivalence relation for non-abelian cohomology coincides with the equivalence relation for Maurer-Cartan elements in the DGLA $(L,[\cdot,\cdot]_C,\overline{\partial})$ (Theorem $\ref{MCtoext2}$).

\section{Classification of non-abelian extensions of Leibniz algebras in terms of non-abelian cohomology}\label{sec:pre}

 {\bf A Leibniz algebra}   is a vector space
$\frkg$, endowed with a linear map
$[\cdot,\cdot]_\frkg:\frkg\otimes\frkg\longrightarrow\frkg$
satisfying
\begin{equation}\label{eq:Leibniz}
~[x,[y,z]_\frkg]_\frkg=[[x,y]_\frkg,z]_\frkg+[y,[x,z]_\frkg]_\frkg,\quad
\forall~x,y,z\in \frkg.
\end{equation}

This is in fact a left Leibniz algebra. In this paper, we only
consider left Leibniz algebras. Denote by $\Z(\g)$ the left
 center of $\g$, i.e.
\begin{equation}\label{leftcen}
\Z(\g)=\{x\in\g\mid [x,y]_\g=0,~~\forall y\in\g\}.
\end{equation}

{\bf A representation} of the Leibniz algebra
$(\frkg,[\cdot,\cdot]_\frkg)$ is a triple  $(V,l,r)$, where  $V$ is
a vector space equipped with two linear maps
$l:\g\longrightarrow\gl(V)$ and $r:\g\longrightarrow\gl(V)$ such
that the following equalities hold:
\begin{equation}\label{condition of rep}
l_{[x,y]_\frkg}=[l_{x},l_{y}],\quad
r_{[x,y]_\frkg}=[l_{x},r_{y}],\quad r_{y}\circ l_{x}=-r_{y}\circ
r_{x}, ~~~ \, ~~~\, ~~~  \forall x,y\in\frkg.
\end{equation}
  The
Leibniz cohomology of $\g$ with coefficients in $V$ is the
cohomology of the cochain complex
$C^k(\g,V)=\Hom(\otimes^{k+1}\g,V), (k\geq0)$ with the
coboundary operator $\partial:C^{k-1}(\g,V)\longrightarrow
C^{k}(\g,V)$ defined by
\begin{eqnarray}
\nonumber\partial
c(x_1,\dots,x_{k+1})&=&\sum_{i=1}^k(-1)^{i+1}l_{x_i}(c(x_1,\dots,\widehat{x_i},\dots,x_{k+1}))
+(-1)^{k+1}r_{x_{k+1}}(c(x_1,\dots,x_k))\\
\label{formulapartial}&&+\sum_{1\leq i<j\leq
k+1}(-1)^ic(x_1,\dots,\widehat{x_i},\dots,x_{j-1},[x_i,x_j]_\g,x_{j+1},\dots,x_{k+1}).
\end{eqnarray}

We  denote by $\Der^L(\g)$ and $\Der^R(\g)$ the set of left derivations and the set of right derivations of $\g$ respectively:
\begin{eqnarray*}
\Der^L(\g)&=&\{D\in \gl(\g)\vert D[x,y]_{\g}=[Dx,y]_{\g}+[x,Dy]_{\g},\forall x,y\in \g\},\\
 \Der^R(\g)&=&\{D\in \gl(\g)\vert D[x,y]_{\g}=[x,Dy]_{\g}-[y,Dx]_{\g},\forall x,y\in \g\}.
 \end{eqnarray*}
It is easy to see that for all $x\in\g$, $\ad^L_x:\g\longrightarrow\g$, which is given by $\ad^L_x(y)=[x,y]_{\g}$, is a left derivation; $\ad^R_x:\g\longrightarrow\g$, which is given by $\ad^R_x(y)=[y,x]_{\g}$, is a right derivation.

\begin{lem}\label{LtoCC}
If $x\in\Z(\g)$,    then for all $D^L\in\Der^L(\g)$, we have $D^Lx\in\Z(\g)$.
\end{lem}
\pf It follows from
$$[D^Lx,y]=D^L[x,y]-[x,D^Ly]=0,\quad \forall~x\in\Z(\g).\qed$$

For all $D_1, D_2\in\gl(\g)$,  $ [D_1,D_2]$ denotes their commutator, i.e. $ [D_1,D_2]=D_1D_2-D_2D_1$. By straightforward computations, we have
\begin{lem}\label{oper}
For all $x\in \g$, $D^L, D^L_1,D^L_2 \in \Der^L(\g)$ and $D^R\in \Der^R(\g)$, we have
\begin{itemize}
\item[$\rm(i)$] $[D^L_1,D^L_2]\in \Der^L(\g) ;$
\item[$\rm(ii)$]$[D^L,D^R]\in \Der^R(\g);$
\item[$\rm(iii)$]$[D^L,\ad^L_x]=\ad^L_{D^Lx},~~[D^L,\ad^R_x]=\ad^R_{D^Lx},~~[\ad^L_x,D^R]=\ad^R_{D^Rx}.$
\end{itemize}
\end{lem}

\begin{cor}
  With the above notations, $\Der^L(\g)$ is a Lie algebra, and $\Der^R(\g)$ is a $\Der^L(\g)$-module, where the representation $\rho$ is given by
  $$\rho(D^L)(D^R)=[D^L,D^R],\quad \forall D^L  \in \Der^L(\g), D^R\in \Der^R(\g).$$
\end{cor}

Denote by $\D(\g)=\Der^L(\g)\oplus \Der^R(\g)$.

View $\Der^L(\g)$ as a Leibniz algebra, then $(\rho,0)$ is a representation of $\Der^L(\g)$ on $\Der^R(\g)$. Thus, we have the semidirect product Leibniz algebra.
\begin{lem}
 $(\D(\g),[\cdot,\cdot]_s)$ is a Leibniz algebra, where the bracket operation $[\cdot,\cdot]_s$ on $\D(\g)$ is given by
\begin{equation}
[(D^L_1,D^R_1),(D^L_2,D^R_2)]_s=([D^L_1,D^L_2],[D^L_1,D^R_2]),\quad \forall~D^L_i\in\Der^L(\g),~~D^R_i\in\Der^R(\g).
\end{equation}
\end{lem}

\begin{defi}\begin{itemize}
\item[\rm(1)] Let $\g$, $\h$, $\hat{\g}$ be  Leibniz algebras. A non-abelian extension of Leibniz algebras is a short exact sequence of Leibniz algebras:
$$ 0\longrightarrow\h\stackrel{\id}{\longrightarrow}\hat{\g}\stackrel{\p}\longrightarrow\g\longrightarrow0.$$
We say that $\hat{\g}$ is a non-abelian  extension of $\g$ by $\h$.
\item[\rm(2)] A splitting of $\hat{\g}$ is a linear map $\sigma:\g\rightarrow\widehat{\g}$ such that $\p\circ \sigma=\idd$.
\item[\rm(3)] Two extensions of $\g$ by $\h$,  $\hat{\g}_1$ and $\hat{\g}_2$,   are said to be isomorphic if there exists a Leibniz algebra morphism $\theta:\hat{\g}_2\longrightarrow \hat{\g}_1$ such that we have the following commutative diagram:
\begin{equation}\label{diagram1}
\begin{array}{ccccccccc}
0&\longrightarrow& \h&\stackrel{\id_2}\longrightarrow&\hat{\g}_2&\stackrel{\p_2}\longrightarrow&\g&\longrightarrow&0\\
 &            &\Big\|&       &\theta\Big\downarrow&          &\Big\|& &\\
 0&\longrightarrow&\h&\stackrel{\id_1}\longrightarrow&\hat{\g}_1&\stackrel{\p_1}\longrightarrow&\g&\longrightarrow&0.
 \end{array}\end{equation}
\end{itemize}
\end{defi}

Let $\hat{\g}$ be a non-abelian  extension of $\g$ by $\h$, and  $\sigma:\g\rightarrow\widehat{\g}$ a splitting. Define $\omega:\g\otimes\g\rightarrow\h$,
 $l:\g\longrightarrow\gl(\h)$ and $r:\g\longrightarrow\gl(\h)$ respectively by
\begin{eqnarray}
  \label{eq:str1}\omega(x,y)&=&[\sigma(x),\sigma(y)]_{\hat{\g}}-\sigma[x,y]_{\g},\quad \forall x,y\in\g,\\
 \label{eq:str2} l_x(\beta)&=&[\sigma(x),\beta]_{\hat{\g}},\quad \forall x\in\g,~\beta\in \h,\\
 \label{eq:str31} r_y(\alpha)&=&[\alpha,\sigma(y)]_{\hat{\g}},\quad \forall y\in\g,~\alpha\in \h.
\end{eqnarray}

Given a splitting, we have $\hat{\g}\cong\g\oplus \h$, and
the Leibniz algebra structure on  $\hat{\g}$ can be transferred to  $\g\oplus\h$:
\begin{eqnarray}\label{bracket of nbe}
[x+\alpha,y+\beta]_{(l,r,\omega)}=[x,y]_{\g}+\omega(x,y)+l_x\beta+r_y\alpha+[\alpha,\beta]_{\h}.\label{bracket of Ext}
\end{eqnarray}

\begin{pro}\label{nonabelian extension of LB}
$(\g\oplus\h,[\cdot,\cdot]_{(l,r,\omega)})$ is a Leibniz algebra if and only if $l,r,\omega$ satisfy the following equalities:
\begin{eqnarray}
l_x[\alpha,\beta]_{\h}&=&[l_x\alpha,\beta]_{\h}+[\alpha,l_x\beta]_{\h},\label{l Der}\\
r_x[\alpha,\beta]_{\h}&=&[\alpha,r_x\beta]_{\h}-[\beta,r_x\alpha]_{\h},\label{r Der}\\
{[l_x\alpha+r_x\alpha,\beta]}_{\h}&=&0,\label{Left Cen1}\\
{[l_x,l_y]}-l_{[x,y]_{\g}}&=&\ad^L_{\omega(x,y)},\label{rep1}\\
{[l_x,r_y]}-r_{[x,y]_{\g}}&=&\ad^R_{\omega(x,y)},\label{rep2}\\
r_y(r_x(\alpha)+l_x(\alpha))&=&0,\label{Left Cen2}\\
l_x\omega(y,z)-l_y\omega(x,z)-r_z\omega(x,y)&=&\omega([x,y]_{\g},z)-\omega(x,[y,z]_{\g})+\omega(y,[x,z]_{\g}).\label{cocycle}
\end{eqnarray}
\end{pro}
\pf Assume that $(\g\oplus\h,[\cdot,\cdot]_{(l,r,\omega)})$ is a Leibniz algebra.   By
$$[x,[\alpha,\beta]_{(l,r,\omega)}]_{(l,r,\omega)}=[[x,\alpha]_{(l,r,\omega)},\beta]_{(l,r,\omega)}+[\alpha,[x,\beta]_{(l,r,\omega)}]_{(l,r,\omega)},$$
we deduce that \eqref{l Der} holds. By
$$[\alpha,[\beta,x]_{(l,r,\omega)}]_{(l,r,\omega)}=[[\alpha,\beta]_{(l,r,\omega)},x]_{(l,r,\omega)}+[\beta,[\alpha,x]_{(l,r,\omega)}]_{(l,r,\omega)},$$
we deduce that  \eqref{r Der} holds.
By
$$[x,[y,\alpha]_{(l,r,\omega)}]_{(l,r,\omega)}=[[x,y]_{(l,r,\omega)},\alpha]_{(l,r,\omega)}+[y,[x,\alpha]_{(l,r,\omega)}]_{(l,r,\omega)},$$
we deduce that \eqref{rep1}   holds.
By
$$[x,[\alpha,y]_{(l,r,\omega)}]_{(l,r,\omega)}=[[x,\alpha]_{(l,r,\omega)},y]_{(l,r,\omega)}+[\alpha,[x,y]_{(l,r,\omega)}]_{(l,r,\omega)},$$
we deduce that \eqref{rep2}   holds.
By
$$[x,[y,z]_{(l,r,\omega)}]_{(l,r,\omega)}=[[x,y]_{(l,r,\omega)},z]_{(l,r,\omega)}+[y,[x,z]_{(l,r,\omega)}]_{(l,r,\omega)},$$
we deduce that \eqref{cocycle}  holds.
By
$$[\alpha,[x,\beta]_{(l,r,\omega)}]_{(l,r,\omega)}=[[\alpha,x]_{(l,r,\omega)},\beta]_{(l,r,\omega)}+[x,[\alpha,\beta]_{(l,r,\omega)}]_{(l,r,\omega)},$$
we obtain that
$l_x[\alpha,\beta]_{\h}=[\alpha,l_x\beta]_{\h}-[r_x\alpha,\beta]_{\h}.$
Combining with \eqref{l Der}, we deduce that \eqref{Left Cen1} holds.
By
$$[\alpha,[x,y]_{(l,r,\omega)}]_{(l,r,\omega)}=[[\alpha,x]_{(l,r,\omega)},y]_{(l,r,\omega)}+[x,[\alpha,y]_{(l,r,\omega)}]_{(l,r,\omega)},$$
we obtain that
${l_xr_y}-r_{[x,y]_{\g}}+r_yr_x=\ad^R_{\omega(x,y)}.$
Combining with \eqref{rep2}, we deduce that \eqref{Left Cen2} holds.

Conversely, if \eqref{l Der}-\eqref{cocycle} hold, it is straightforward to see that $(\g\oplus\h,[\cdot,\cdot]_{(l,r,\omega)})$ is a Leibniz algebra. The proof is finished.
\qed\vspace{3mm}

Note that \eqref{l Der} means that $l_x\in \Der^L(\h)$, and \eqref{r Der} means that $r_x\in \Der^R(\h)$.

\begin{defi}\label{2cocyle}
A non-abelian $2$-cocycle on $\g$ with values in $\h$ is a triple $(l,r,\omega)$ of linear maps $l:\g\rightarrow \Der^L(\h),~~r:\g\rightarrow \Der^R(\h)$ and $\omega:\g\otimes\g\rightarrow \h$ satisfying Eqs. \eqref{Left Cen1}-\eqref{cocycle}.
We denote by $Z^2(\g,\h)$ the set of non-abelian $2$-cocycles.

Two $2$-cocycles $(l^1,r^1,\omega^1)$ and $(l^2,r^2,\omega^2)$ are equivalent, \emptycomment{$(l^1,r^1,\omega^1)\stackrel{\varphi}\sim(l^2,r^2,\omega^2)$,} if there exists a linear map $\varphi:\g\rightarrow\h$ such that  for all $x,y\in\g$, the following equalities hold:
\begin{eqnarray}
l^1_x-l^2_x&=&\ad^L_{\varphi(x)},\label{cohomology 1}\\
r^1_x-r^2_x&=&\ad^R_{\varphi(x)},\label{cohomology 2}\\
\omega^1(x,y)-\omega^2(x,y)&=&l^2_x\varphi(y)+r^2_y\varphi(x)+[\varphi(x),\varphi(y)]_{\h}-\varphi([x,y]_{\g})\label{cohomology 3}.
\end{eqnarray}

The second non-abelian cohomology $H^2(\g,\h)$ is the quotient of $Z^2(\g,\h)$ by the above equivalence relation.
\end{defi}

\begin{thm}\label{classification}
Let $(\g,[\cdot,\cdot]_\g)$ and $(\h,[\cdot,\cdot]_\h)$ be Leibniz algebras. Then non-abelian extensions of $\g$ by $\h$ are classified by the second non-abelian cohomology $H^2(\g,\h)$.
\end{thm}
\pf Let $\hat{\g}$ be a non-abelian extension of $\g$ by $\h$. By choosing a splitting $\sigma_1:\g\longrightarrow\hat{\g}$, we obtain a 2-cocycle $(l^1,r^1,\omega^1)$.  First we show that the cohomological class of $(l^1,r^1,\omega^1)$ does not depend on the choice of splittings. In fact,
let $\sigma_1$ and $\sigma_2$ be two different splittings. Define $\varphi:\g\longrightarrow\h$ by $\varphi(x):=\sigma_1(x)-\sigma_2(x)$.
 By $l^j_x(\beta)=[\sigma_j(x),\beta]_{\hat{\g}}$ and $r^j_y(\alpha)=[\alpha,\sigma_j(y)]_{\hat{\g}},~~j=1,2$, it is obvious that    $(\ref{cohomology 1})$ and $(\ref{cohomology 2})$ hold.
 Furthermore, we have
\begin{eqnarray*}
\omega^1(x,y)&=&[\sigma_1(x),\sigma_1(y)]_{\hat{\g}}-\sigma_1[x,y]_{\g}\\
&=&[\sigma_2(x)+\varphi(x),\sigma_2(y)+\varphi(y)]_{\hat{\g}}-\sigma_2[x,y]_{\g}-\varphi([x,y]_{\g})\\
&=&[\sigma_2(x),\sigma_2(y)]_{\hat{\g}}+[\sigma_2(x),\varphi(y)]_{\hat{\g}}+[\varphi(x),\sigma_2(y)]_{\hat{\g}}+[\varphi(x),\varphi(y)]_{\h}-\sigma_2[x,y]_{\g}-\varphi[x,y]_{\g}\\
&=&\sigma_2[x,y]_{\g}+\omega^2(x,y)+l^2_x\varphi(y)+r^2_y\varphi(x)+[\varphi(x),\varphi(y)]_{\h}-\sigma_2[x,y]_{\g}-\varphi([x,y]_{\g})\\
&=&\omega^2(x,y)+l^2_x\varphi(y)+r^2_y\varphi(x)+[\varphi(x),\varphi(y)]_{\h}-\varphi([x,y]_{\g}),
\end{eqnarray*}
which implies that $(\ref{cohomology 3})$ holds. Thus, $(l^1,r^1,\omega^1)$ and $(l^2,r^2,\omega^2)$  are in the same cohomological class.

Now we go on to prove that isomorphic extensions give rise to the same element in $H^2(\g,\h)$. Assume that $\hat{\g_1}$ and $\hat{\g_2}$ are isomorphic
  extensions of $\g$ by $\h$, and $\theta:\hat{\g_2}\longrightarrow\hat{\g_1}$ is a Leibniz algebra morphism such that we have the commutative diagram
  \eqref{diagram1}. Assume that $\sigma_2:\g\longrightarrow\hat{\g_2}$ and $\sigma_1:\g\longrightarrow \hat{\g_1}$ are two splittings.
Define  $\sigma'_2:\g\longrightarrow\hat{\g_2}$ by $\sigma'_2=\theta^{-1}\circ \sigma_1$. Since $\p_1\circ\theta=\p_2$, we have
$$\p_2\circ \sigma'_2(x)=(\p_1\circ\theta)\circ\theta^{-1}\circ\sigma_1(x)=x.$$
Thus, $\sigma'_2$ is a splitting of $\hat{\g_2}$.  Define $\varphi_\theta:\g\longrightarrow\h$ by $\varphi_{\theta}(x)=\sigma'_2(x)-\sigma_2(x)$.
Since $\theta^{-1}:\g_1\rightarrow\g_2$ is a morphism of Leibniz algebras,  we have
\begin{eqnarray*}
\theta^{-1}{[\sigma_1(x),\alpha]}_{\hat{\g}_1}={[\theta^{-1}\sigma_1(x),\alpha]}_{\hat{\g}_2}
={[\varphi_{\theta}(x)+\sigma_2(x),\alpha]}_{\hat{\g}_2}
=l^2_x\alpha+\ad^L_{\varphi_{\theta}(x)}\alpha.
\end{eqnarray*}
On the other hand, since $\theta|_\h=\idd$, we have $$\theta^{-1}[\sigma_1(x),\alpha]_{\hat{\g}_1}=[\sigma_1(x),\alpha]_{\hat{\g}_1}=l^1_x\alpha.$$
Therefore, we have
\begin{equation}\label{eq:t1}
l^1_x-l^2_x=\ad^L_{\varphi_{\theta}(x)}.
\end{equation}
Similarly, we have
\begin{equation}\label{eq:t2}
r^1_x-r^2_x=\ad^R_{\varphi_{\theta}(x)}.
\end{equation}
For all $x,y\in\g$, on one hand,  we have
\begin{eqnarray*}
\theta^{-1}{[\sigma_1(x),\sigma_1(y)]}_{\hat{\g}_1}&=&{[\theta^{-1}\sigma_1(x),\theta^{-1}\sigma_1(y)]}_{\hat{\g}_2}\\
&=&{[\varphi_{\theta}(x)+\sigma_2(x),\varphi_{\theta}(y)+\sigma_2(y)]}_{\hat{\g}_2}\\
&=&[\varphi_{\theta}(x),\varphi_{\theta}(y)]_{\h}+l^2_x\varphi_{\theta}(y)+r^2_y\varphi_{\theta}(x)+[\sigma_2(x),\sigma_2(y)]_{\hat{\g}_2}\\
&=&[\varphi_{\theta}(x),\varphi_{\theta}(y)]_{\h}+l^2_x\varphi_{\theta}(y)+r^2_y\varphi_{\theta}(x)+\sigma_2[x,y]_{\g}+\omega^2(x,y).
\end{eqnarray*}
On the other hand, we have
\begin{eqnarray*}
\theta^{-1}{[\sigma_1(x),\sigma_1(y)]}_{\hat{\g}_1}&=&\theta^{-1}(\sigma_1[x,y]_{\g}+\omega^1(x,y))=
\varphi_{\theta}[x,y]_{\g}+\sigma_2[x,y]_{\g}+\omega^1(x,y).
\end{eqnarray*}
In conclusion, we have
\begin{equation}\label{eq:t3}
\omega^1(x,y)-\omega^2(x,y)=l^2_x\varphi_{\theta}(y)+r^2_y\varphi_{\theta}(x)+[\varphi_{\theta}(x),\varphi_{\theta}(y)]_{\h}-\varphi_{\theta}[x,y]_{\g}.
\end{equation}
By \eqref{eq:t1}-\eqref{eq:t3}, we deduce that $(l^1,r^1,\omega^1)$ and $(l^2,r^2,\omega^2)$ are in the same cohomological class.

Conversely, given two 2-cocycles $(l^1,r^1,\omega^1)$ and  $(l^2,r^2,\omega^2)$ representing the same cohomological class, i.e. there exists $\varphi:\g\longrightarrow\h$ such that \eqref{cohomology 1}-\eqref{cohomology 3} holds, we show that the corresponding extensions $(\g\oplus\h,[\cdot,\cdot]_{(l^1,r^1,\omega^1)})$ and $(\g\oplus\h,[\cdot,\cdot]_{(l^2,r^2,\omega^2)})$ are isomorphic.
 Define $\theta:\g\oplus\h\longrightarrow\g\oplus\h$ by $$\theta(x+\alpha)=x-\varphi(x)+\alpha.$$
 Then, we can deduce that $\theta$ is a Leibniz algebra morphism such that the following diagram commutes:
 \begin{equation*}
\begin{array}{ccccccccc}
0&\longrightarrow&\h&\stackrel{\id_2}\longrightarrow&\g\oplus\h_{(l^2,r^2,\omega^2)}&\stackrel{\p_2}\longrightarrow&\g&\longrightarrow&0\\
 &            &\Big\|&       &\theta\Big\downarrow&          &\Big\|& &\\
 0&\longrightarrow&\h&\stackrel{\id_1}\longrightarrow&\g\oplus\h_{(l^1,r^1,\omega^1)}&\stackrel{\p_1}\longrightarrow&\g&\longrightarrow&0.
 \end{array}\end{equation*}
We omit details. This finishes the proof.\qed

\section{Another approach to non-abelian extensions of Leibniz algebras }\label{Leib20ext}

The notion of strongly homotopy (sh) Leibniz algebras, or
$Lod_\infty$-algebras was introduced in \cite{livernet}, see also $\cite{ammardefiLeibnizalgebra,UchinoshL}$ for more details. In \cite{Leibniz2al}, the authors introduced the notion of a Leibniz 2-algebra, which is the categorification of a Leibniz algebra, and prove that the category of Leibniz 2-algebras and the category of 2-term $Lod_\infty$-algebras are equivalent.


\emptycomment{
\begin{defi}$\cite{ammardefiLeibnizalgebra}$
A sh Leibniz algebra is a graded  vector space $L=L_0\oplus
L_1\oplus\cdots$ equipped with a system $\{l_k|~1\leq k<\infty\}$ of
linear maps $l_k:\wedge^kL\longrightarrow L$ with degree
$\deg(l_k)=k-2$, where the exterior powers are interpreted in the
graded sense and the following relation  is satisfied:
\begin{eqnarray*}
&&\sum_{i+j=Const}\sum_{k\geq
j}\sum_{\sigma}(-1)^{(k+1-j)(j-1)}(-1)^{j(|x_{\sigma(1)}|+\cdots+|x_{\sigma(k-j)}|)}\sum_{\sigma}\sgn(\sigma)\Ksgn(\sigma)\\
&&l_i(x_{\sigma(1)},\dots,x_{\sigma(k-j)},l_j(x_{\sigma(k+1-j)},\dots,x_{\sigma(k)}),x_{\sigma(k+1)},\dots,x_{\sigma(i+j-1)})=0,
\end{eqnarray*}
where the summation is taken over all $(k-j,j-1)$-unshuffles and
``$\Ksgn(\sigma)$'' is the Koszul sign for a permutation $\sigma\in
S_k$, i.e. $$ x_1\wedge x_2\wedge\cdots\wedge
x_k=\Ksgn(\sigma)x_{\sigma(1)}\wedge x_{\sigma(2)}\wedge\cdots\wedge
x_{\sigma(k)}.
$$
\end{defi}

In particular, if we concentrate on the 2-term case, we can give
explicit formulas for 2-term sh Leibniz algebras as follows:
}

\begin{defi}\label{defi:2leibniz}
  A   Leibniz $2$-algebra $\huaV$ consists of the following data:
\begin{itemize}
\item[$\bullet$] a complex of vector spaces $\huaV:V_1\stackrel{\dM}{\longrightarrow}V_0,$

\item[$\bullet$] bilinear maps $l_2:V_i\times V_j\longrightarrow
V_{i+j}$, where $0\leq i+j\leq 1$,

\item[$\bullet$] a  trilinear map $l_3:V_0\times V_0\times V_0\longrightarrow
V_1$,
   \end{itemize}
   such that for any $w,x,y,z\in V_0$ and $m,n\in V_1$, the following equalities are satisfied:
\begin{itemize}
\item[$\rm(a)$] $\dM l_2(x,m)=l_2(x,\dM m),$
\item[$\rm(b)$]$\dM l_2(m,x)=l_2(\dM m,x),$
\item[$\rm(c)$]$l_2(\dM m,n)=l_2(m,\dM n),$
\item[$\rm(d)$]$\dM l_3(x,y,z)=l_2(x,l_2(y,z))-l_2(l_2(x,y),z)-l_2(y,l_2(x,z)),$
\item[$\rm(e_1)$]$ l_3(x,y,\dM m)=l_2(x,l_2(y,m))-l_2(l_2(x,y),m)-l_2(y,l_2(x,m)),$
\item[$\rm(e_2)$]$ l_3(x,\dM m,y)=l_2(x,l_2(m,y))-l_2(l_2(x,m),y)-l_2(m,l_2(x,y)),$
\item[$\rm(e_3)$]$ l_3(\dM m,x,y)=l_2(m,l_2(x,y))-l_2(l_2(m,x),y)-l_2(x,l_2(m,y)),$
\item[$\rm(f)$] the Jacobiator identity:\begin{eqnarray*}
&&l_2(w,l_3(x,y,z))-l_2(x,l_3(w,y,z))+l_2(y,l_3(w,x,z))+l_2(l_3(w,x,y),z)\\
&&-l_3(l_2(w,x),y,z)-l_3(x,l_2(w,y),z)-l_3(x,y,l_2(w,z))\\
&&+l_3(w,l_2(x,y),z)+l_3(w,y,l_2(x,z))-l_3(w,x,l_2(y,z))=0.\end{eqnarray*}
   \end{itemize}
\end{defi}
We usually denote a   Leibniz 2-algebra by
$(V_1,V_0,\dM,l_2,l_3)$, or simply by
$\huaV$. A Leibniz 2-algebra is a {\bf strict Leibniz 2-algebra} if $l_3=0.$ If $l_2$ is skewsymmetric, we obtain {\bf Lie 2-algebras}, see \cite{baez:2algebras} for more details.

\begin{lem}\label{lem:rep1}
For a  Leibniz $2$-algebra
$(V_1,V_0,\dM,l_2,l_3)$, we have
\begin{equation}\label{eqn:rep1}
l_2(l_2(x,m),y)+l_2(l_2(m,x),y)=0,\quad\forall~x,y\in V_0,m\in V_1.
\end{equation}
\end{lem}

\begin{defi}\label{defi:Leibniz morphism}
 Let $\huaV$ and $\huaV^\prime$ be  Leibniz $2$-algebras, a morphism $\frkf$
 from $\huaV$ to $\huaV^\prime$ consists of
\begin{itemize}
  \item[$\bullet$] linear maps $f_0:V_0\longrightarrow V_0^\prime$
  and $f_1:V_1\longrightarrow V_1^\prime$ commuting with the
  differential, i.e.
  $$
f_0\circ \dM=\dM^\prime\circ f_1;
  $$
  \item[$\bullet$] a bilinear map $f_2:V_0\times V_0\longrightarrow
  V_1^\prime$,
\end{itemize}
  such that  for all $x,y,z\in L_0,~m\in
L_1$, we have
\begin{equation}\label{eqn:DGLA morphism c 1}\left\{\begin{array}{rll}
l_2^\prime(f_0(x),f_0(y))-f_0l_2(x,y)&=&\dM^\prime
f_2(x,y),\\
l_2^\prime(f_0(x),f_1(m))-f_1l_2(x,m)&=&f_2(x,\dM
m),\\
l_2^\prime(f_1(m),f_0(x))-f_1l_2(m,x)&=&f_2(\dM
m,x),\end{array}\right.
\end{equation}
and
\begin{eqnarray}
\nonumber&&f_1(l_3(x,y,z))+l_2^\prime(f_0(x),f_2(y,z))-l_2^\prime(f_0(y),f_2(x,z))-l_2^\prime(f_2(x,y),f_0(z))\\
\label{eqn:DGLA morphism c 2}
&&-f_2(l_2(x,y),z)+f_2(x,l_2(y,z))-f_2(y,l_2(x,z))-l_3^\prime(f_0(x),f_0(y),f_0(z))=0.
\end{eqnarray}
\end{defi}

Let $(\h,[\cdot,\cdot]_\h)$ be a Leibniz algebra. Define $\Delta(\h)\subset \Z(\h)$ by
\begin{equation}
\Delta(\h)=\mathrm {span}\{[\alpha,\alpha]_\g\mid\alpha \in \h\}.
\end{equation}
Define  $\Pi(\h)\subset\D(\h)$ by
\begin{equation}
\Pi(\h)=\{(D^L,D^R)\in\D(\h)\mid D^L\alpha+D^R\alpha\in\Delta(\h),~\forall \alpha\in\h\}.
\end{equation}

\begin{lem}\label{lem:DR}
  For any $\beta\in\Delta(\h)$ and $D^R\in\Der^R(\h)$, we have $D^R\beta=0$.
\end{lem}
\pf  For all $\alpha\in\h$ and $D^R\in\Der^R(\h)$, we have  $D^R[\alpha,\alpha]_\h=0$, which implies that for all $\beta\in\Delta(\h)$, $D^R\beta=0$.\qed

\begin{pro}\label{pro:subalgebra}
$\Pi(\h)$ is a sub-Leibniz algebra of $(\D(\h),[\cdot,\cdot]_s)$.
\end{pro}
\pf By Lemma \ref{lem:DR}, for any $\beta\in\Delta(\h)$ and $(D^L,D^R)\in\Pi(\h)$, $D^L\beta\in\Delta(\h)$.
Thus, for all $(D^L_1,D^R_1),(D^L_2,D^R_2)\in \Pi(\h)$, we have
\begin{eqnarray*}
[(D^L_1,D^R_1),(D^L_2,D^R_2)]_s\alpha&=&[D^L_1,D^L_2]\alpha+[D^L_1,D^R_2]\alpha\\
&=&D^L_1(D_2^L\alpha+D_2^R\alpha) -(D_2^LD^L_1\alpha+D_2^RD^L_1\alpha)\in\Delta(\h),
\end{eqnarray*}
 which implies that $[(D^L_1,D^R_1),(D^L_2,D^R_2)]_s\in\Pi(\h)$.
Thus, $\Pi(\h)$ is a sub-Leibniz algebra.\qed\vspace{3mm}

It is obvious that for all $\alpha\in\h$, $(\ad^L_\alpha,\ad^R_\alpha)\in\Pi(\h)$.

Consider the complex $\h\stackrel{(\ad^L,\ad^R)}{\longrightarrow}\Pi(\h)$, and define $l_2$ by
\begin{equation}\label{eq:l2}
\left\{\begin{array}{rlll}
l_2((D^L_1,D^R_1),(D^L_2,D^R_2))&=&[(D^L_1,D^R_1),(D^L_2,D^R_2)]_s,&\quad\mbox{for $(D^L_i,D^R_i) \in \Pi(\h)$},\\
l_2((D^L,D^R),\alpha)&=&D^L\alpha,&\quad\mbox{for $(D^L,D^R) \in \Pi(\h),~\alpha\in\h$},\\
l_2(\alpha,(D^L,D^R))&=&D^R\alpha,&\quad\mbox{for $(D^L,D^R) \in \Pi(\h),~\alpha\in\h$}.
\end{array}\right.
\end{equation}

\begin{pro}\label{pro:Leibniz2}
With the above notations, $(\h,\Pi(\h),(\ad^L,\ad^R),l_2)$ is a strict Leibniz $2$-algebra, where $l_2$ is given by \eqref{eq:l2}.
\end{pro}
\pf
By Lemma $\ref{oper}$, we have
\begin{eqnarray*}
(\ad^L,\ad^R) l_2((D^L,D^R),\alpha)&=&(\ad^L,\ad^R) (D^L\alpha)=\ad^L_{D^L\alpha}+\ad^R_{D^L\alpha}\\
&=&[D^L,\ad^L_\alpha]+[D^L,\ad^R_\alpha]\\
&=&[(D^L,D^R),(\ad^L,\ad^R) \alpha]_s=l_2((D^L,D^R),(\ad^L,\ad^R)\alpha),
\end{eqnarray*}
which implies that Condition (a) in Definition \ref{defi:2leibniz} holds.

Similarly, we have
\begin{eqnarray*}
\dM l_2(\alpha,D)-l_2(\dM\alpha,D)&=&\ad^L_{D^R\alpha}+\ad^R_{D^R\alpha}-([\ad^L_\alpha,D^L]+[\ad^L_\alpha,D^R])\\
&=&\ad^L_{D^L\alpha+D^R\alpha} =0,
\end{eqnarray*}
which implies that Condition (b) in Definition \ref{defi:2leibniz} holds.

Condition (c) follows from
$$
l_2((\ad^L_\alpha,\ad^R_\alpha),\beta)=l_2(\alpha,(\ad^L_\beta,\ad^R_\beta))=[\alpha,\beta]_\h.
$$

By Proposition \ref{pro:subalgebra}, Condition (d) holds naturally.

For all $(D^L_1,D^R_1),(D^L_2,D^R_2)\in \Pi(\h)$,  we have
\begin{eqnarray*}
l_2(\alpha,l_2((D^L_1,D^R_1),(D^L_2,D^R_2)))-l_2(l_2(\alpha,(D^L_1,D^R_1)),(D^L_2,D^R_2))-l_2((D^L_1,D^R_1),l_2(\alpha,(D^L_2,D^R_2)))\\
=[D^L_1,D^R_2]\alpha-D^R_2D^R_1\alpha-D^L_1D^R_2\alpha
=-D^R_2(D_1^L\alpha+D_1^R\alpha)=0,
\end{eqnarray*}
which implies that Condition (${\rm e_1}$) in Definition \ref{defi:2leibniz} holds.

Similarly, Conditions (${\rm e_2}$) and (${\rm e_3}$) also hold. Since $l_3=0$, Condition (f) holds naturally. The proof is finished.
\qed\vspace{3mm}

It is well-known that
when $\h$ is a Lie algebra, there is a strict Lie 2-algebra $(\h,\Der(\h),\ad,l_2)$, where $\Der(\h)$ is the set of derivations of $\h$, and $l_2$ is given by
\begin{equation}
\left\{\begin{array}{rlll}
l_2(D_1,D_2)&=&[D_1,D_2],&\quad\mbox{for $D_1,D_2\in \Der(\h)$},\\
l_2(D,\alpha)&=&-l_2(\alpha,D)=D\alpha,&\quad\mbox{for $D\in \Der(\h),~\alpha\in\h$}.
\end{array}\right.
\end{equation}
When $\h$ is a Lie algebra, it is obvious that $\Delta(\h)=0$. Therefore, $\Pi(\h)$ reduces to
$$
\Pi(\h)=\{(D,-D)|D\in\Der(\h)\}.
$$
Furthermore, the strict Leibniz 2-algebra given in Proposition \ref{pro:Leibniz2} reduces to a Lie 2-algebra. Define $f_0:\Der(\h)\longrightarrow\Pi(\h)$ by
$$f_0:D\longmapsto (D,-D),\quad \forall D\in \Der(\h).$$ About the relation between the two Lie 2-algebras, we have

\begin{pro}
Let $\h$ be a Lie algebra, then $(f_0,f_1=\idd,f_2=0)$ is an isomorphism from the Lie $2$-algebra $(\h,\Der(\h),\ad,l_2)$   to $(\h,\Pi(\h),(\ad^L,\ad^R),l_2)$.
\end{pro}
\pf It is obvious that $f_0$ and $f_1$ commute with the differential. Furthermore, we have
\begin{eqnarray*}
  f_0l_2(D_1,D_2)&=&([D_1,D_2],-[D_1,D_2])=[(D_1,-D_1),(D_2,-D_2)]_s=l_2(f_0(D_1),f_0(D_2)),\\
  f_1l_2(D,\alpha)&=&f_1(D\alpha)=D\alpha=l_2(f_0(D),f_1(\alpha)).
\end{eqnarray*}
 Thus, $(f_0,f_1=\idd,f_2=0)$ is an isomorphism between Lie 2-algebras. \qed

 \begin{rmk}
   From the above result, we see that the Leibniz $2$-algebra associated to a Leibniz algebra given in Proposition \ref{pro:Leibniz2} is a natural generalization of the Lie $2$-algebra associated to a Lie algebra. However, we find that we can construct another Leibniz $2$-algebra associated to a Leibniz algebra, which plays an important role when we study non-abelian extensions of Leibniz algebras.
 \end{rmk}

 Define a subspace $\Xi(\h)\subset \D(\h)$ by
 \begin{equation}
 \Xi(\h)=\{(D^L,D^R)\in\D (\h)\mid D^L\alpha+D^R\alpha\in \Z(\h),~\mbox{for all $\alpha\in\h$,~and}~ D^R\Z(\h)=0\}.
 \end{equation}
 It is obvious that for all $\alpha\in\h$, $(\ad^L_\alpha,\ad^R_\alpha)\in\Xi(\h)$.

\begin{pro}
$\Xi(\h)$ is a sub-Leibniz algebra of  $(\D (\h),[\cdot,\cdot]_s)$.
\end{pro}
\pf
For all $(D^L_1,D^R_1),~(D^L_2,D^R_2)\in\Xi(\h)$ and $\alpha,\beta\in\h$, by Lemma $\ref{LtoCC}$, we have
\begin{eqnarray*}
[[(D^L_1,D^R_1),(D^L_2,D^R_2)]_s\alpha,\beta]_\h&=&[[D^L_1,D^L_2]\alpha+[D^L_1,D^R_2]\alpha,\beta]_\h\\
&=&[D^L_1(D_2^L\alpha+D_2^R\alpha),\beta]_\h-[D_2^LD^L_1\alpha+D_2^RD^L_1\alpha,\beta]_\h\\
&=&0.
\end{eqnarray*}
Furthermore, by Lemma $\ref{LtoCC}$ and the fact that $D^R_2\Z(\h)=0$, for all $\alpha\in\Z(\h)$, we have
$$[D^L_1,D^R_2]\alpha=D^L_1D^R_2\alpha-D^R_2D^L_1\alpha=0.$$
Thus, $[(D^L_1,D^R_1),(D^L_2,D^R_2)]_s\alpha\in\Xi(\h),$ which implies that  $\Xi(\h)$ is a sub-Leibniz algebra.\qed

\begin{cor}
$\Xi(\h)$ is a right ideal of $(\D(\h),[\cdot,\cdot]_s)$.
\end{cor}

Similar as Proposition \ref{pro:Leibniz2}, we have

\begin{pro}\label{Leibniz 2 algebra}
With the above notations, $(\h,\Xi(\h),(\ad^L,\ad^R),l_2)$ is a strict Leibniz $2$-algebra, where $l_2$ is given by \eqref{eq:l2} for all $(D^L_i,D^R_i) \in \Xi(\h),~(D^L,D^R)\in\Xi(\h).$
\end{pro}

Define $\Pro_L:\D(\h)\longrightarrow\Der^L(\h)$ and $\Pro_R:\D(\h)\longrightarrow\Der^R(\h)$   by
$$\Pro_L(D^L,D^R)=D^L,~~~\Pro_R(D^L,D^R)=D^R.$$

\begin{thm}\label{LB2toext1}
Let $\hat{\g}=(\g\oplus\h,[\cdot,\cdot]_{(l,r,\omega)})$ be a non-abelian extension of $\g$ by $\h$ with $\Z(\h)=\Z(\hat{\g})\cap \h$. Then, $(l,r,\omega)$ give rise to a Leibniz $2$-algebra morphism $(f_0,f_1,f_2)$ from   $\g$ to the Leibniz $2$-algebra $(\h,\Xi(\h),(\ad^L,\ad^R),l_2)$,
where $f_0,~~f_1$ and $f_2$ are given by
$$f_0(x)=(l_x,r_x),~~f_1=0,~~f_2(x,y)=\omega(x,y).$$
Conversely, for any morphism $(f_0,f_1,f_2)$ from   $\g$ to   $(\h,\Xi(\h),(\ad^L,\ad^R),l_2)$, there is a non-abelian  extension  $\hat{\g}=(\g\oplus\h,[\cdot,\cdot]_{(l,r,\omega)})$ of $\g$ by $\h$, with $\Z(\h)=\Z(\hat{\g})\cap \h$, where   $l,r,\omega$ are given by
$$l_x=\Pro_Lf_0(x),~~~r_x=\Pro_Rf_0(x),~~~\omega(x,y)=f_2(x,y).$$
\end{thm}
\pf Let $\hat{\g}=(\g\oplus\h,[\cdot,\cdot]_{(l,r,\omega)})$ be a non-abelian extension of $\g$ by $\h$ with $\Z(\h)=\Z(\hat{\g})\cap \h$.
Then, $l,r,\omega$ satisfy Eqs. \eqref{l Der}-\eqref{cocycle}. By $(\ref{l Der})$ and $(\ref{r Der})$, we have $l_x\in\Der^L(\g)$ and $r_x\in\Der^R(\g)$.
By $(\ref{Left Cen1})$ and $(\ref{Left Cen2})$, we have $l_x(\alpha)+r_x(\alpha)\in \Z(\hat{\g})\cap\h\subset \Z(\h)$, which implies that $f_0(x)(\alpha)\in\Z(\h)$, for all $\alpha\in \h$. Furthermore, it is obvious that    $r_x(\Z(\hat{\g})\cap\h)=0$. By the condition
$\Z(\h)=\Z(\hat{\g})\cap \h$, we get $r_x(\Z(\h))=0$. Thus, we deduce that $f_0(x)\in\Xi(\h)$. Furthermore, by $(\ref{rep1})$ and $(\ref{rep2})$, we  get
\begin{eqnarray*}
l_2(f_0(x),f_0(y))-f_0([x,y]_\g)&=&([l_x,l_y]-l_{[x,y]_{\g}},[l_x,r_y]-r_{[x,y]_{\g}})=(\ad^L_{\omega(x,y)},\ad^R_{\omega(x,y)})\\
&=&(\ad^L,\ad^R)f_2(x,y).
\end{eqnarray*}
At last, by $(\ref{cocycle})$, we deduce that \eqref{eqn:DGLA morphism c 2} holds. Thus, $(f_0,f_1,f_2)$ is a morphism from $\g$ to the Leibniz 2-algebra $(\h,\Xi(\h),(\ad^L,\ad^R),l_2)$.

Conversely,  by the fact that $(f_0,f_1,f_2)$ is a morphism from $\g$ to  $(\h,\Xi(\h),(\ad^L,\ad^R),l_2)$, it is straightforward to deduce that Eqs.
$(\ref{l Der})$-$(\ref{cocycle})$ hold. Furthermore, for all $\alpha\in\Z(\h)$, by the fact that $r_y\Z(\h)=0$, we have
 $$[\alpha,y+\beta]_{(l,r,\omega)}=r_y\alpha+[\alpha,\beta]_\h=0,$$
 i.e. $\alpha\in\Z(\hat{\g})\cap \h$. Thus, we have $\Z(\h)=\Z(\hat{\g})\cap \h$.
 The proof is finished.\qed
\begin{rmk}
  We can also prove that isomorphic extensions correspond to equivalent Leibniz $2$-algebra morphisms similar as Theorem \ref{classification}. We omit details.
\end{rmk}

By the results in \cite{Class of 4}, we can give some examples of non-abelian extensions of Leibniz algebras.

\begin{ex}{\rm
Let $\h=\{e_1,e_2,e_3\}$ be the Leibniz algebra with the non-zero brackets $[e_1,e_2]=e_3$, and $\g=\{x\}$  the trivial Leibniz algebra. Let $\hat{\g}=\{e_1,e_2,e_3,x\}$ be the Leibniz algebra with the non-zero brackets
\begin{eqnarray*}
 [e_1,e_2]=e_3,~~[x,e_1]=e_1,~~[x,e_2]=\gamma e_2, ~~
 [x,e_3]=(1+\gamma)e_3,~~[e_1,x]=-e_1,~~\gamma \in\mathbb R.
\end{eqnarray*}
Then one has
$$\Z(\hat{\g})=\Z(\h)=\{e_2,e_3\}.$$}
\end{ex}

\begin{ex}{\rm
Let $\h=\{e_1,e_2,e_3\}$ be the Leibniz algebra with the non-zero brackets $[e_1,e_1]=e_3$, and $\g=\{x\}$  the trivial Leibniz algebra. Let $\hat{\g}=\{e_1,e_2,e_3,x\}$ be the Leibniz algebra with the non-zero brackets
\begin{eqnarray*}
[e_1,e_1]=e_3,~~[x,e_2]=e_2,~~[e_2,x]=-e_2.
\end{eqnarray*}
Then we have
$$\Z(\h)=\{e_2,e_3\},~~~\Z(\hat{\g})=\{e_3\}.$$
Note that $\Z(\hat{\g})\subsetneqq \Z(\h)$.}
\end{ex}

\begin{ex}{\rm
Let $\h=\{e_1,e_2,e_3\}$ be the Leibniz algebra with the non-zero brackets $[e_1,e_1]=e_3$, and $\g=\{x\}$ the trivial Leibniz algebra. Let $\hat{\g}_1=\{e_1,e_2,e_3,x\}$ be the Leibniz algebra with the non-zero brackets
\begin{eqnarray*}
 [e_1,e_1]=e_3,[x,e_1]=e_1+e_2,~~[x,e_2]=e_2,~~
 [x,e_3]=2e_3,~~[e_1,x]=-e_1-e_2,~~[e_2,x]=-e_2.
\end{eqnarray*}
Let $\hat{\g}_2=\{e_1,e_2,e_3,x\}$ be the Leibniz algebra with the non-zero brackets
$$[e_1,e_1]=e_3,~~[x,e_2]=e_2.$$

Then one has
$\Z(\h)=\{e_2,e_3\}$. But $\Z(\hat{\g}_2)=\{e_2,e_3\}$ and $\Z(\hat{\g}_1)=\{e_3\}$. They are not equal.
Therefore, this two extensions are not equivalent.}
\end{ex}

\emptycomment{

Furthermore, we have
$$r^1_xl^2_xe_2=r^x_1e_2=-e_2.$$

 From this example, we can see it is impossible to find a unified Leibniz $2$-algebra to describe all extensions of a Leibniz algebra $\hat{\g}$ of $\g$ by $\h$. This also indicates that why we can't discus a general case of extension.

}

\section{Non-abelian extensions of Leibniz algebras and Maurer-Cartan elements}\label{dgla0ext}
In this section, we generalize the results in \cite{nonabelin cohomology of Lie} to the case of Leibniz algebras. Let $(\g,[\cdot,\cdot]_\g)$ be a Leibniz algebra. The graded vector space $C(\g,\g)=\oplus_kC^k(\g,\g)$ equipped with the graded commutator bracket
\begin{equation}\label{eq:bracket}
[P,Q]_{C}=P\circ Q+(-1)^{pq+1}Q\circ P,~~~\forall P\in C^{p}(\g,\g),~~Q\in C^{q}(\g,\g)
\end{equation}
is a graded Lie algebra, where $C^k(\g,\g)=\Hom(\otimes^{k+1}\g,\g)$ and $P\circ Q\in C^{p+q}(\g,\g)$ is defined by
\begin{eqnarray*}
P\circ Q(x_1,\cdots,x_{p+q+1})&=&\sum^p_{k=0}(-1)^{kq}(\sum_{\sigma\in sh(k,q)}\sgn(\sigma)P(x_{\sigma(1)},\cdots,x_{\sigma(k)},\\
&& Q(x_{\sigma(k+1)},\cdots,x_{\sigma(k+q)},x_{k+q+1}),x_{k+q+2},\cdots,x_{p+q+1})).
\end{eqnarray*}
Furthermore, $(C(\g,\g),[\cdot,\cdot]_C,\overline{\partial})$ is a DGLA, where $\overline{\partial }$ is given by $\overline{\partial}P=(-1)^p\partial P$ for all $P\in C^{p}(\g,\g)$, and $\partial $ is the coboundary operator of $\g$ with the coefficients in the adjoint representation $(\ad^L,\ad^R)$. See \cite{dgLa1,dgLa2}  for more details.
\begin{rmk}
The coboundary operator $\partial$ associated to the adjoint representation of the Leibniz algebra $\g$ can be written as $\partial P=(-1)^p[\mu_\g,P]_{C}$, for all $P\in C^{p}(\g,\g)$,   where $\mu_{\g}\in C^{1}(\g,\g) $ is the Leibniz algebra structure on $\g$, i.e.   $\mu_\g(x,y)=[x,y]_{\g}$. Thus, we have $\overline{\partial}P=[\mu_\g,P]_{C}.$
\end{rmk}

The set $MC(L)$ of {\bf Maurer-Cartan elements} of a DGLA $(L,[\cdot,\cdot],\dM)$ is defined by
$$MC(L)\triangleq \{P\in L_1\mid\dM P+\frac{1}{2}[P,P]=0\}.$$

Let $(\g,[\cdot,\cdot]_\g)$ and $(\h,[\cdot,\cdot]_\h)$ be two Leibniz algebras. Let $\g\oplus\h$ be the Leibniz algebra direct sum of   $\g$ and $\h$, where the bracket is defined by  $[x+\alpha,y+\beta]=[x,y]_\g+[\alpha,\beta]_\h$. Then there is DGLA  $(C(\g\oplus\h,\g\oplus\h), [\cdot,\cdot]_{C},\overline{\partial})$, where $\partial $ is the coboundary operator for the  Leibniz algebra $\g\oplus\h$ with the coefficients in the adjoint representation. It is straightforward to see that $(C(\g\oplus\h, \h), [\cdot,\cdot]_{C},\overline{\partial})$ is a sub-DGLA. Define $C^k_>(\g\oplus\h,\h)\subset C^k(\g\oplus\h,\h)$ by
$$C^k(\g\oplus\h,\h)={C}^k_>(\g\oplus\h,\h)\oplus C^k(\h,\h).$$
  Denote by $C_>(\g\oplus\h,\h)=\oplus_kC^k_>(\g\oplus\h,\h)$, which is a graded vector space.


\begin{lem}\label{lem:dgla}
  With the above notations, we have $\partial (C^k_>(\g\oplus\h,\h))\subset C^{k+1}_>(\g\oplus\h,\h)$, and $({C}_>(\g\oplus\h,\h),[\cdot,\cdot]_{C},\overline{\partial})$ is a sub-DGLA of $(C(\g\oplus\h,\g\oplus\h), [\cdot,\cdot]_{C},\overline{\partial})$.    Furthermore, its degree $0$ part $C^0_>(\g\oplus\h,\h)=\Hom(\g,\h)$ is abelian.
\end{lem}
We   denote by $(L=\oplus_k L_k,[\cdot,\cdot]_C,\overline{\partial})$ this differential graded Lie algebra,   where $L_k=C^k_>(\g\oplus\h,\h).$

\pf By the fact that $\overline{\partial}= [\mu_{\g}+\mu_{\h},\cdot]_{C}$, we can verify that $\overline{\partial }(C^k_>(\g\oplus\h,\h))\subset C^{k+1}_>(\g\oplus\h,\h)$. It is straightforward to see that
${C}_>(\g\oplus\h,\h)$ is closed under the bracket operation ${[\cdot,\cdot]}_{C}$. Obviously, $C^0_>(\g\oplus\h,\h)=\Hom(\g,\h)$ is abelian.  \qed
\vspace{3mm}

We have $L_1=\Hom(\otimes ^2\g,\h)\oplus\Hom(\g\otimes \h,\h)\oplus\Hom(\h\otimes \g,\h)$. For all $\varphi_1,\varphi_2\in L_0,~~P\in L_1$,   we have
\begin{equation}\label{eq:001}
[[P,\varphi_1]_C,\varphi_2]_C=0.
\end{equation}
\begin{defi}
Let $(L,[\cdot,\cdot]_C,\overline{\partial})$ be the DGLA given in Lemma \ref{lem:dgla}. Two elements $c,c'\in L_1$ are said to be equivalent
if there is an element $\varphi\in L_0$ such that
\begin{equation}\label{eq:equirelation}
c'= c+[c,\varphi]_C+\overline{\partial}\varphi+\frac{1}{2}{[\overline{\partial}\varphi,\varphi]}_{C}.
\end{equation}
\end{defi}

\begin{rmk}
By \eqref{eq:001} and the fact that $L_0$ is abelian, we can deduce that \eqref{eq:equirelation} is indeed an equivalence relation. See \cite{nonabelin cohomology of Lie} for more details about this equivalence  relation  for Lie algebras.
\end{rmk}

 Now consider $(\g\oplus\h,[\cdot,\cdot]_{(l,r,\omega)})$, where $[\cdot,\cdot]_{(l,r,\omega)}$ is given by \eqref{bracket of nbe} for any $l,r:\g\longrightarrow\gl(\h)$ and $\omega\in\Hom(\otimes^2\g,\h)$. View $l,r:\g\longrightarrow\gl(\h)$ as elements in $\Hom(\g\otimes \h,\h)$ and $\Hom(\h\otimes \g,\h)$  respectively via
\begin{equation}\label{eq:tensorequivalent}
 l(x,\alpha)=l_x(\alpha),\quad r(\alpha,x)=r_x(\alpha).
\end{equation}

\begin{pro}\label{pro:JMC}
The following two statements are equivalent:
\begin{itemize}
  \item[\rm(a)] $(\g\oplus\h,[\cdot,\cdot]_{(l,r,\omega)})$ is a Leibniz algebra, which is a non-abelian extension of $\g$ by $\h$;

  \item[\rm(b)] $l+r+\omega$ is a Maurer-Cartan element, i.e. $ \omega+l+r\in MC(L).$
\end{itemize}

\end{pro}
\pf
By  Proposition \ref{nonabelian extension of LB}, $(\g\oplus\h,[\cdot,\cdot]_{(l,r,\omega)})$ is a Leibniz algebra if and only if Eqs. \eqref{l Der}-\eqref{cocycle} hold.

If $c=\omega+l+r$ is a Maurer-Cartan element, we have
$$(\overline{\partial} c+\frac{1}{2}[c,c]_C)(e_1,e_2,e_3)=0, \quad \forall e_i=x_i+\alpha_i\in\g\oplus \h.$$
By straightforward computations, we have
\begin{eqnarray*}
 \partial c(e_1,e_2,e_3)&=&[e_1,c(e_2,e_3)] -[e_2,c(e_1,e_3)] -[c(e_1,e_2),e_3] \\
 &&-c([e_1,e_2],e_3)-c(e_2,[e_1,e_3])+c(e_1,[e_2,e_3]),
 \end{eqnarray*}
and
\begin{eqnarray*}
{[e_1,c(e_2,e_3)]} &=&{[\alpha_1,\omega(x_2,x_3)]}_{\h}+{[\alpha_1,l_{x_2}\alpha_3]}_{\h}+{[\alpha_1,r_{x_3}\alpha_2]}_{\h},\\
{[e_2,c(e_1,e_3)]} &=&{[\alpha_2,\omega(x_1,x_3)]}_{\h}+{[\alpha_2,l_{x_1}\alpha_3]}_{\h}+{[\alpha_2,r_{x_3}\alpha_1]}_{\h},\\
{[c(e_1,e_2),e_3]}&=&{[\omega(x_1,x_2),\alpha_3]}_{\h}+{[l_{x_1}\alpha_2,\alpha_3]}_{\h}+{[r_{x_2}\alpha_1,\alpha_3]}_{\h},\\
c([e_1,e_2],e_3)&=&\omega([x_1,x_2]_\g,x_3)+l_{[x_1,x_2]_\g}\alpha_3+r_{x_3}[\alpha_1,\alpha_2]_\h,\\
c(e_2,[e_1,e_3])&=&\omega(x_2,[x_1,x_3]_\g)+l_{x_2}[\alpha_1,\alpha_3]_\h+r_{[x_1,x_3]_\g}\alpha_2,\\
c(e_1,[e_2,e_3])&=&\omega(x_1,[x_2,x_3]_\g)+l_{x_1}[\alpha_2,\alpha_3]_\h+r_{[x_2,x_3]_\g}\alpha_1.
\end{eqnarray*}
Furthermore, by \eqref{eq:bracket}, we have
\begin{eqnarray*}
\half [c,c]_C=(c\circ c)(e_1,e_2,e_3)&=&c(c(e_1,e_2),e_3)-c(e_1,c(e_2,e_3))+c(e_2,c(e_1,e_3))\\
&=&r_{x_3}\omega(x_1,x_2)+r_{x_3}l_{x_1}\alpha_2+r_{x_3}r_{x_2}\alpha_1-l_{x_1}\omega(x_2,x_3)\\
&&-l_{x_1}l_{x_2}\alpha_3-l_{x_1}r_{x_3}\alpha_2+l_{x_2}\omega(x_1,x_3)+l_{x_2}l_{x_1}\alpha_3+l_{x_2}r_{x_3}\alpha_1.
\end{eqnarray*}
Thus, $c=\omega+l+r$ is a Maurer-Cartan element if and only if \eqref{l Der}-\eqref{cocycle} hold. The proof is finished.\qed

\begin{cor}\label{MCtoext1}
Let $\g$ and $\h$ be two Leibniz algebras, then there is a one-to-one correspondence between  $Z^2(\g,\h)$ and Maurer-Cartan elements in the DGLA $(L,[\cdot,\cdot]_C,\overline{\partial}).$
\end{cor}

\emptycomment{
\pf By Proposition $\ref{nonabelian extension of LB}$, the Jacobi Identity  is equivalent to conditions of the Definition $\ref{2cocyle}$ of non-abelian $2$-cocycles, which means the Jacobi Identity  is equivalent to a non-abelian $2$-cocycle. However, by Proposition $\ref{JMC}$, the Jacobi Identity is also equivalent to one Maurer-Cartan element. In a conclusion, one non-abelian $2$-cocycle is equivalent to one Maurer-Cartan element.\qed
}

We can further prove that

\begin{thm}\label{MCtoext2}
Let $\g$ and $\h$ be two Leibniz algebras, the equivalence in $Z^2(\g,\h)$ coincides with the equivalence in $MC(L)$.
\end{thm}
\pf
Two elements $c=\omega+r+l$ and $c'=\omega'+r'+l'$ in $MC(L)$ are equivalent if there exists $\varphi\in \Hom(\g,\h)$ such that
$$c'= c+[c,\varphi]_C+\partial \varphi+{[\partial \varphi,\varphi]}_{C}.$$
For all $e_i=x_i+\alpha_i\in\g\oplus\h$, we have
\begin{eqnarray*}
[c,\varphi]_C(e_1,e_2)&=&[\omega+r+l,\varphi]_{C}(e_1,e_2)
= l_{x_1}\varphi(x_2)+r_{x_2}\varphi(x_1),\\
\partial\varphi(e_1,e_2)&=&{[\alpha_1,\varphi(x_2)]}_{\h}+{[\varphi(x_1),\alpha_2]}_{\h}-\varphi([x_1,x_2]_\g),\\
{[\partial\varphi,\varphi]}_{C}(e_1,e_2)&=&2{[\varphi(x_1),\varphi(x_2)]}_{\h}.
\end{eqnarray*}
Therefore, we have
$$c'(e_1,e_2)=c(e_1,e_2)+ l_{x_1}\varphi(x_2)+r_{x_2}\varphi(x_1)+\ad^R_{\varphi(x_2)}\alpha_1+\ad^L_{\varphi(x_1)}\alpha_2-\varphi([x_1,x_2]_\g)+{[\varphi(x_1),\varphi(x_2)]}_{\h},$$
which implies that
\begin{eqnarray*}
  l'(x_1,\alpha_2)-l(x_1,\alpha_2)&=&\ad^L_{\varphi(x_1)}\alpha_2,\\
   r'(\alpha_1,x_2)-l(\alpha_1,x_2)&=&\ad^R_{\varphi(x_2)}\alpha_1,\\
   \omega'(x_1,x_2)-\omega(x_1,x_2)&=&l_{x_1}\varphi(x_2)+r_{x_2}\varphi(x_1)-\varphi([x_1,x_2]_\g)+{[\varphi(x_1),\varphi(x_2)]}_{\h}.
\end{eqnarray*}

By \eqref{eq:tensorequivalent}, they are equivalent to
Eqs. (\ref{cohomology 1})-(\ref{cohomology 3}), which are the equivalence relation in $Z^2(\g,\h)$. This finishes the proof.\qed

\end{document}